# A profit allocation mechanism for customer-to-manufacturer platform in e-commercial business

Bo DAI[1], Fenfen LI[2*]
[1]*School of Management, Hunan University of Technology and Business, Changsha, 410205, China*
[2]*School of Public Administration, Hunan University of Finance and Economics, Changsha, 410205, China*
[1]*dr_bo_dai@163.com*, [2]*lifenfen@hufe.edu.cn*

**Abstract:**
A centralized collaboration problem of customer-to-manufacturer (C2M) platform in e-commercial business is studied in this paper, where an e-commercial company and multiple manufacturers come to an agreement for constituting a collaborative alliance to satisfy the product orders of multiple customers. The problem has two important issues which are the optimal reallocation of customer orders among the manufacturers to maximize a total profit of the alliance and a fair allocation of the profit among the manufacturers so as to assure the stability of the alliance and the sustainability of C2M platform. A profit allocation mechanism is proposed which is based on core allocation concept and considers the contribution of each manufacturer in fulfilling orders. The effectiveness of the mechanism is evaluated with numerical experiments on sixty instances generated based on real data of Alibaba.

## 1 Introduction

Nowadays, with the changing of customers' consumption pattern, the fierce competition, the increase of labour cost and the increase of raw material cost, how to produce high-quality products that meet the variable and customized needs of customers in a lower cost has become a very crucial issue for e-commercial companies and manufacturers. This has brought great challenges to the traditional supply chain where there are several intermediate levels between customers and manufacturers, such as wholesalers and retailers. As we know, e-commerce business is booming in the past ten years. In recent years, a new e-commercial business model called customer to manufacturer (C2M for short) has appeared. C2M connects to the customer in one end, and connects to the manufacturer in the other end. C2M subverts the traditional retail thinking, where manufacturing is driven by customer demands through the e-commerce platform. C2M is a new ecommerce variant that links users and the factory directly to eliminate all unnecessary channels (Luo et al., 2019).

There are some evident advantages to C2M. Firstly, C2M bring cost savings for both customers and manufacturers. C2M realizes the direct connection from customers to the manufacturers and removes the middlemen and possible price increase links. Thus, customers can obtain good quality products with civilian price. Manufacturers can reduce a large amount of inventory cost based on the pre-sale orders from customers which naturally reduce the cost of the customers to purchase the products. Secondly, C2M is a customer demand driven mode that realizes private customized production as customers can participate in the production process. Manufacturers can achieve a fast respond to the changing of customers' consumption pattern.

Because of the above advantages, C2M business model has gradually attracted attentions in recent years. Several large manufacturing enterprises of China, such as the Red Collar Group (Hu et al., 2016) and Haier (Sun et al., 2015), had utilised the C2M model to attract more end consumers. In addition, Zhang et al. (2019) presented computational framework for designing and optimizing custom compression casts/braces. They proposed a customer-to-manufacturer design model that started from a 3D scanned human model represented by mesh and ended with the 3D printed cast/brace.

However, there are still some difficulties that restrict the implement of C2M. It takes time for manufacturers to produce products after receiving pre-sale orders, which lengthens the process cycle of products to customers. The customers' buying interest will be weakened in a long cycle. In addition, if the order quantity of the C2M platform is not enough for the manufacturer, it is not profitable for him. In order to ensure the on time delivery of small orders, the manufacturer may be at a loss and reluctant to stay in the C2M platform. Therefore, the horizontal and vertical cooperation among the members of C2M platform are an effective solution to solve the above difficulties. Through vertical cooperation with big e-commercial company, the manufacturers can obtain larger orders. Through horizontal cooperation among manufacturers, they can fulfil more production orders in less process cycle. Therefore, cooperation is a promising way for the implement of C2M.

The most important issue for the cooperation of C2M platform is to allocate the total profit among all manufacturers reasonably in a collaborative alliance (grand coalition) so that the alliance is beneficial to all participants and stable and the C2M platform is sustainable. This issue can be regarded as a profit allocation problem, which is usually studied in the framework of cooperative game theory. Many methods have been proposed to deal with this issue, such as the Core concept (Gillies, 1959), Shapley value (Shapley, 1953), and Kalai-Smorodinsky (KS for short) solution (Kalai, 1977), etc. The Core concept is the unique allocation that is stable. If an allocation is not in the core, some players in the game may be leave from the alliance so as to seek a larger profit. Both Shapley value and KS solution cannot guarantee that their allocations are always stable. Furthermore, profit allocation problem have been intensively studied in many areas, such as in collaborative logistics where multiple shippers and carriers cooperate each other in transportation planning (Krajewska et al. 2008; Audy et al., 2011; Dai and Chen, 2012, 2015; Gansterer and Hartl, 2018), in supply chain (Jaber et al., 2006; Frascatore and Mahmoodi, 2011; Yao and Ran, 2019), in virtual enterprises (Chen et al., 2007; Deng and Liu, 2009), etc. But the studies of profit allocation problem in C2M platform are rarely.

In literature, most studies about C2M generally focused on strategic planning and empirical analysis (Gallino et al., 2014 ; Hu et al., 2016). Wang et al. (2018) discussed the motivation and evolution mechanism of business model innovation. To deal with the practical problems the enterprises faces under the "internet +" situation, they proposed a C2M (customer to maker) model that balanced the clients needs and the products and services of the enterprises, where a linkage platform mechanism was constructed between enterprises and customers. But they just studied the C2M business model innovation in the field of strategic management not the profit allocation problem. Liu et al. (2018) studtied the information sharing between the manufacturer and the distribution enterprise in a C2M e-commerce supply chain. They showed that a well-designed revenue reallocation mechanism can improve the efficiency of information sharing, achieve the rationality and effectiveness of supply chain profit allocation and maximize the value of information sharing in the C2M supply chain. But they considered the simple profit function and the game model of Stackelberg not a cooperative game model with integrating core concept.

In this paper, we study a centralized collaboration problem of customer-to-manufacturer platform in e-commercial business, where an e-commercial company and multiple manufacturers set up a collaborative alliance. In collaboration, the orders fulfilled by each manufacturer are determined by the e-commercial company through solving a centralized production planning problem that is formulated as an integer programming model. After collaboration, the issue of fairly allocating the post-collaboration total profit of the manufacturers is solved by proposing a $\gamma$-approximate profit allocation mechanism. The mechanism considers the contribution of each manufacturer in fulfilling orders, and it is evaluated by numerical experiments on sixty instances.

The novelty and contribution of our mechanism can be found in three aspects:

1) To our best knowledge, the profit allocation problem happened in C2M platform is rarely investigated in previous studies.

2) A stable allocation can be provided by our proposed mechanism that integrates the core concept, which can assure the sustainability of C2M platform.

3) The proposed mechanism considers the contribution of each manufacturer in fulfilling orders, which reflects the fairness of the profit allocation.

The rest of this paper is organized as follows. Section 2 presents respectively the studied problem and introduces a centralized collaboration framework for the problem. Section 3 proposes a profit allocation mechanism and relevant models. Section 4 evaluates the performances of the mechanism by numerical experiments. Section 5 concludes this paper with comments for future research.

## 2 Problem description

In the studied customer to manufacturer (C2M) business platform, an e-commercial company and multiple manufacturers come to an agreement for constituting a collaborative alliance to satisfy the product orders of multiple customers. The collaboration is under a centralized collaboration framework. The e-commercial company is in charge of receiving, organizing and fulfiling the product orders of customers, which acts as a coordinator in the alliance. The manufacturers are responsible for fabricating the required products from the coordinator. The collaborative alliance is given in Fig. 1.

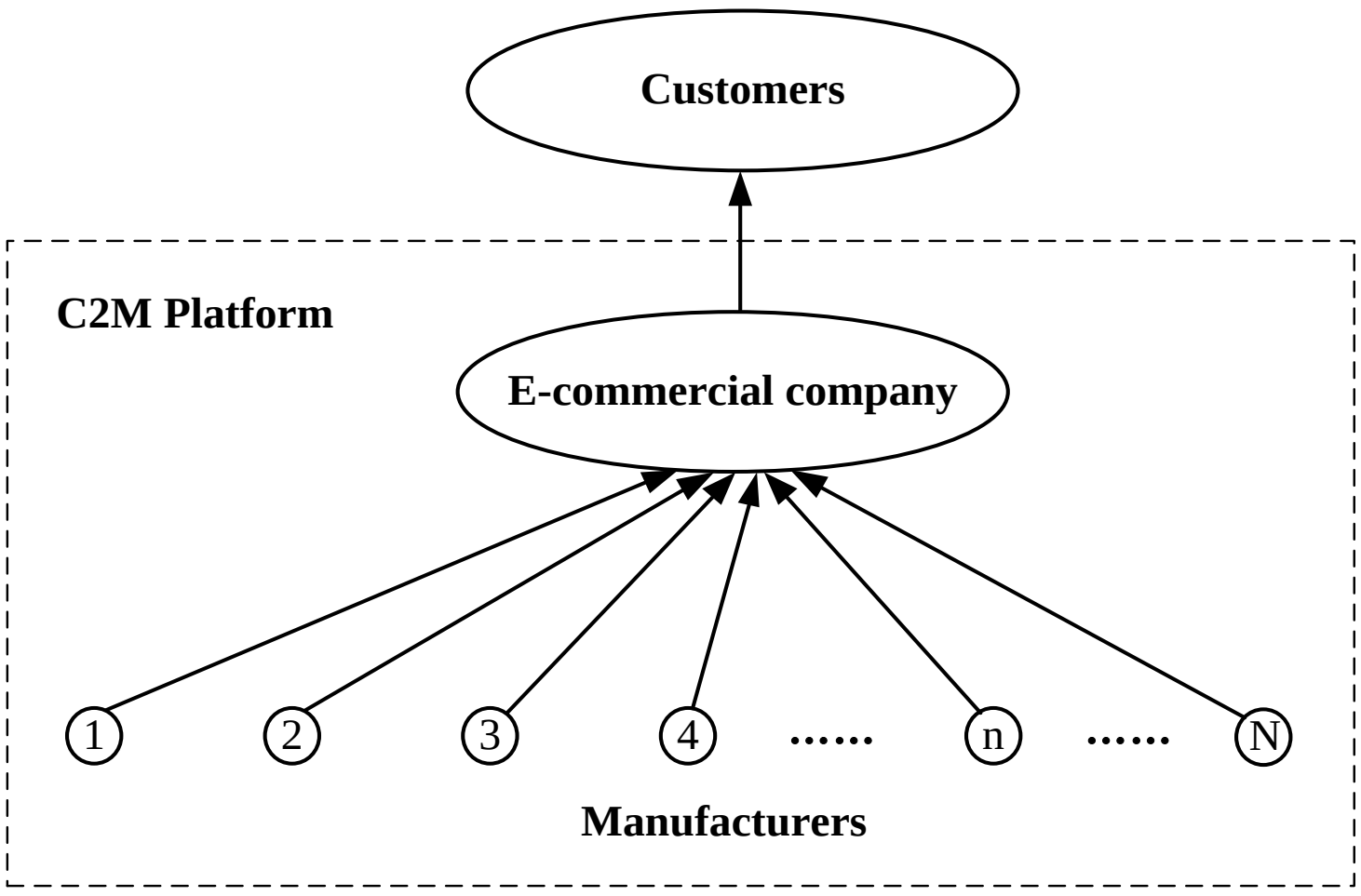

Fig. 1. The collaborative alliance in customer to manufacturer business

Initially, the coordinator has acquired certain product orders from the customers through the C2M platform. Each order (demand) is consist of an order quantity and an order delivery time (days) of each product, where a sale price of each unit product paid by a customer to the coordinator for obtaining the product and a ask price of each unit product paid by the coordinator to the manufacturer for producing the product are also given.

In detail, the framework consists of the following steps.

1) Before collaboration, the coordinator and the manufacturers have signed a binding cooperation agreement. Each manufacturer submits its manufacturing cost of fabricating each product and its production capacity of each product per day to the coordinator.

2) During collaboration, a centralized collaborative production planning is made by the coordinator to maximize a total profit of the alliance, in which the product orders are allocated among manufacturers by the coordinator. Thus, the post-collaboration total profit of the alliance can be generated. Here, shortage cost of per unit each product is estimated by the coordinator.

3) After collaboration, based on the collaborative production plan, the total profit of the alliance is allocated among the manufacturers by using a profit allocation mechanism.

The proposed centralized collaboration framework for C2M business is realistic in the practical cases. This has been explored by some famous e-commercial companies of China. For example, Alibaba has established C2M business unit in its Taobao business group at the end of 2019. Then Taobao has built a new online platform with C2M customized products as the core supply in 2020. In addition, other Chines e-commercial companies such as JD, NetEase have also established the similar platforms. These e-commerce giants already have certain market and customers, so it is easier to establish such platforms for them. Thus, manufacturers are more willing to cooperate with them. The possibility of obtaining profit through collaboration may attract manufacturers to accept a centrally coordinated solution if the profit allocation among the manufacturers is fair. That is also the original intention of the profit allocation mechanism design for C2M platform in this paper.

From above discussion, two important issues need to be addressed for the cooperation in C2M platform. The first issue is a collaborative production planning whose objective is to find a production plan for each manufacturer that satisfy the order quantity and an order delivery time of each order, the production capacity constraints of each manufacturer so that the total profit of the alliance is maximized. The second issue is to allocate the total profit among all manufacturers in a fair way so that the post-collaboration profit of each manufacturer is no less than its pre-collaboration profit and all manufacturers are willing to stay in the alliance. In this paper, we mainly focus on the second issue.

For simplicity of formulation, some assumptions are made for the studied problem. Firstly, all customer orders are allocated among manufacturers by the coordinator. Secondly, all customer orders must be fulfilled under the constraints of the production capacities of the manufacturers, which mean that not all orders can be fulfilled. Thirdly, each order can be fulfilled by multiple manufacturers. Fourthly, different manufacturers may have different manufacturing costs and different production capacities. Finally, the transportation cost and time of each order are not considered as they are organized by the coordinator and completed by other logistics companies.

Note that the profit allocation mechanism proposed in this paper is applicable to other kind of C2M platforms without these assumptions. The introduction of the assumptions is only for formulating the centralized collaborative production planning model in next section.

## 3 Profit allocation mechanisms for C2M platform

With the collaboration, the collaborative alliance can generate a profit called post-collaboration total profit defined by the optimal objective value of the centralized production planning model for C2M platform in this section. This profit will be fairly allocated among all manufacturers in the alliance by a profit allocation mechanism (model) designed in this section. The notation to be used in the two models is given as follows.

*Indices*

$i = 1,\ldots,I$, product index, $I$ is the number of products

$n,m = 1,\ldots,N$, manufacturer index, $N$ is the number of manufacturers

*Parameters*

$MC_{in}$, manufacturing cost of product $i$ of manufacturer $n$

$PC_{in}$, production capacity of product $i$ of manufacturer $n$ per day

$SC_i$, shortage cost per unit product $i$

$SP_i$, ask price per unit product $i$

$OQ_i$, order quantity of product $i$ from customers

$OT_i$, order delivery time (days) of product $i$ from customers

$TP_N$, post-collaboration total profit of C2M platform

$TP_{SN}$, post-collaboration total profit of subset of manufacturers, $SN \in N$

*Variables*

$s_i$, shortage quantity of product $i$ which is an integer

$q_{in}$, quantity of product $i$ allocated to manufacturer $n$ by the coordinator in collaboration which is an integer

$p_n$, post-collaboration profit of manufacturer $n$

$\gamma$, a positive real number, $\gamma \in [0,1]$

*3.1 Centralized production planning model for C2M platform*

With above notation, a centralized production planning model *CPP* for C2M platform can be formulated as the following linear integer programming model.

Model *CPP*:

$$Max\left\{\sum_{n=1}^{N}\sum_{i=1}^{I} q_{in}\cdot(SP_i - MC_{in}) - \sum_{i=1}^{I} s_i \cdot SC_i\right\} \tag{1}$$

Subject to:

$$\sum_{n=1}^{N} q_{in} + s_i = OQ_i, i = 1,\ldots,I \tag{2}$$

$$q_{in} \le PC_{in}\cdot OT_i, i = 1,\ldots,I, n = 1,\ldots,N \tag{3}$$

$$s_i, q_{in} \ge 0, s_i, q_{in} \in Z, i = 1,\ldots,I, n = 1,\ldots,N \tag{4}$$

Objective (1) aims to maximize the post-collaboration total profit of the alliance. Constraints (2) assure that the order quantity of each product is equal to the sum of the allocated quantity of each product to each manufacturer and the shortage quantity of each product. Constraints (3) indicate the production capacity constraint of each product for each manufacturer. Constraints (4) define the value ranges of all variables.

*3.2 Profit allocation model for C2M platform*

After solving model *CPP*, the post-collaboration profit of the alliance is generated, which needs to be fairly allocated among manufacturers. Thus, a post-collaboration profit allocation problem (*PAP*) happened in C2M platform. The cooperative game theory is widely used method to solve the problem in reference. Many solution concepts have been proposed, such as Core concept (Gillies, 1959), Shapley value (Shapley, 1953), Kalai-Smorodinsky solution (Kalai, 1977). The core is the set of imputations (Kalai, 2008) that no coalition has a value larger than the sum of its members' payoffs. If an imputation is in the core, no coalition has incentive to leave the grand coalition (alliance) and acquire a larger payoff. Thus, the core allocations are stable. Sadly, neither Shapley value allocation nor Kalai-Smorodinsky allocation can always stay in the core. Therefore, the two kinds of allocations may be not stable. Besides, the core maybe empty so that the core allocation is not existent. To deal with the case of empty core, Caprara and Letchford (2010) proposed a γ-budget balanced allocation method.

Inspiration by above methods, we design a $\gamma$-approximate core guaranteed profit allocation mechanism in this paper. Our mechanism can generate a stable profit allocation when the core of the problem is empty. If the core is empty, our mechanism can generate a $\gamma$-approximate non empty core in which γ is a value in range [0, 1].

Model *PAP*:

$$Max\{\gamma\} \tag{5}$$

Subject to:

$$p_n \geq 0, n = 1, \dots, N \tag{6}$$

$$\sum_{n=1}^{N} p_n = TP_N \tag{7}$$

$$\sum_{n \in SN} p_n \geq \gamma \cdot TP_{SN}, \forall TP_{SN} \in N \tag{8}$$

$$If \ \sum_{i=1}^{K} \overline{q_{in}} \cdot SP_i \geq \sum_{i=1}^{K} \overline{q_{im}} \cdot SP_i, \ then \quad p_n \geq p_m, n, m = 1, \dots, N, n \neq m \tag{9}$$

$$If \ \sum_{i=1}^{K} \overline{q_{in}} \cdot SP_i \leq \sum_{i=1}^{K} \overline{q_{im}} \cdot SP_i, \ then \quad p_n \leq p_m, n, m = 1, \dots, N, n \neq m \tag{10}$$

$$\gamma \in (0,1], p_n \in R, n = 1, \dots, N \tag{11}$$

Objective (5) aims to maximize $\gamma$ and assure that the profit allocation is always existent. Constraints (6) ensure the individual rationality of the profit allocation. Constraint (7) makes sure that the allocation is efficient or budget-balanced. Constraints (8) guarantee the stability of the grand coalition (alliance). Constraints (6), (7) and (8) together ensure that the allocation is in the core. Constraints (9) and (10) assure that a manufacturer who fulfils more valuable production orders will be allocated more or at least the same post-collaboration profit compared with other manufacturers. Here, $\overline{q_{in}}$ denote the quantity of product *i* allocated to manufacturer *n* by the coordinator obtained through solving model *CPP*. Constraints (11) define the value ranges of all variables.

After solving model *PAP*, a $\gamma$-approximate non empty core allocation can be found. Thus, our mechanism is relative fair and keeps the stability of the alliance for C2M platform in some sense. Note that, pre-collaboration profit of manufacturer $n$ is zero as it will not gain profit if he dosen't participate in collaboration.

## 4 Numeric experiments

In this section, we conduct numerical experiments to evaluate the performances of our proposed profit allocation mechanism. Six sets of instances are generated based on real data of Alibaba. The ask price of each product , the shortage cost of each product, the order (demand) quantitiy of each product are set in the same way with our previous paper (Dai et al., 2019, 2020).

For the instances generated, the number of products is set to 1 ($I$ = 1), 5 ($I$ = 5) or 10 ($I$ = 10), the number of manufacturers is set to 5 ($I$ = 5) or 10 ($N$ = 10). By combing 3 possible numbers of manufacturers, 2 possible numbers of products, we generated 6 sets of instances as indicated in the following table, where each set contains 10 instances. Then each set of instances is indicated by its number of manufacturers and its number of products. For example, the instance set 1x5 contains ten instances with 1 product and 5 manufacturers, and the instance set 1x10 contains ten instances with 1 product and 10 manufacturers.

Table 1: Six sets of instances

| Instance set | Number of products | Number of manufacturers |
|---|---|---|
| 1 | 1 | 5 |
| 2 | 1 | 10 |
| 3 | 5 | 5 |
| 4 | 5 | 10 |
| 5 | 10 | 5 |
| 6 | 10 | 10 |

All models involved in the numerical experiments are solved by using CPLEX 12.9 on a personal PC with i7-8565U CPU and 16GB RAM. As the computation time of each instance is quite short and usually in several seconds, it is not recorded in these tables. The computational results of all instances are given in Table 1 to Table 8. In the tables, column "Post-collaboration total profit" represents the total profit of model *CPP* with objective (1) found by CPLEX MIP solver. Column "The value of $\gamma$" denotes the value of model *PAP* with objective (2) found by CPLEX MIP solver. Column "Product shortage" denotes that whether there is any product shortage in the solution of model *PAP*, where "Y" means product shortage and "N" means no product shortage.

Table 2 Results of the instances set 1 (1x5)

| **Instance** | **1** | **2** | **3** | **4** | **5** | **6** | **7** | **8** | **9** | **10** |
|---|---|---|---|---|---|---|---|---|---|---|
| Post-collaboration total profit | 98954 | 42153 | 138713 | 52874 | 16569 | 18532 | 11432 | 144294 | 59810 | 92058 |
| The value of $\gamma$ | 1 | 1 | 1 | 1 | 1 | 1 | 1 | 1 | 1 | 1 |
| Product shortage | Y | Y | Y | Y | Y | Y | Y | Y | Y | Y |

Table 3 Results of the instances set 2 (1x10)

| **Instance** | **1** | **2** | **3** | **4** | **5** | **6** | **7** | **8** | **9** | **10** |
|---|---|---|---|---|---|---|---|---|---|---|
| Post-collaboration total profit | 218730 | 99960 | 272705 | 199590 | 92937 | 251420 | 30900 | 58020 | 149595 | 185790 |
| The value of $\gamma$ | 0.77 | 0.74 | 0.76 | 0.81 | 1 | 0.76 | 0.9 | 0.93 | 0.79 | 0.76 |
| Product shortage | N | N | N | N | Y | N | N | N | N | N |

Table 4 Results of the instances set 3 (5x5)

| **Instance** | **1** | **2** | **3** | **4** | **5** | **6** | **7** | **8** | **9** | **10** |
|---|---|---|---|---|---|---|---|---|---|---|
| Post-collaboration total profit | 164085 | 250815 | 346496 | 207749 | 335169 | 288368 | 156792 | 337251 | 311377 | 502110 |
| The value of $\gamma$ | 1 | 1 | 1 | 1 | 1 | 1 | 1 | 1 | 1 | 1 |
| Product shortage | Y | Y | Y | Y | Y | Y | Y | Y | Y | Y |

Table 5 Results of the instances set 4 (5x10)

| **Instance** | **1** | **2** | **3** | **4** | **5** | **6** | **7** | **8** | **9** | **10** |
|---|---|---|---|---|---|---|---|---|---|---|
| Post-collaboration total profit | 574109 | 785705 | 741715 | 625500 | 902842 | 654620 | 506242 | 892727 | 902565 | 1033810 |
| The value of $\gamma$ | 0.98 | 0.95 | 0.88 | 0.89 | 0.93 | 0.86 | 0.91 | 0.95 | 0.99 | 0.77 |
| Product shortage | Y | N | N | N | Y | N | Y | Y | N | N |

Table 6 Results of the instances set 5 (10x5)

| **Instance** | **1** | **2** | **3** | **4** | **5** | **6** | **7** | **8** | **9** | **10** |
|---|---|---|---|---|---|---|---|---|---|---|
| Post-collaboration total profit | 739620 | 522703 | 646087 | 610641 | 593944 | 589579 | 668083 | 438621 | 658213 | 416059 |
| The value of $\gamma$ | 1 | 1 | 1 | 1 | 1 | 1 | 1 | 1 | 1 | 1 |
| Product shortage | Y | Y | Y | Y | Y | Y | Y | Y | Y | Y |

Table 7 Results of the instances set 6 (10x10)

| **Instance** | **1** | **2** | **3** | **4** | **5** | **6** | **7** | **8** | **9** | **10** |
|---|---|---|---|---|---|---|---|---|---|---|
| Post-collaboration total profit | 2275625 | 1297824 | 1513747 | 1405410 | 1290160 | 1670663 | 1945640 | 1441654 | 2193560 | 937045 |
| The value of $\gamma$ | 0.97 | 0.85 | 0.87 | 0.86 | 0.84 | 0.92 | 0.93 | 0.96 | 0.95 | 0.88 |
| Product shortage | N | Y | Y | N | N | Y | N | Y | Y | N |

As shown in Tables 2 to 7, we can observe the following results.

1) If the number of products is given, the post-collaboration total profit of the C2M platform increases notablely along with the increase of the number of manufacturers, and the product shortage decreases or even disappears in most instances. Therefore, the participation of more manufacturers in the C2M platform are profitable.

2) If the number of manufacturers is given, the post-collaboration total profit of the C2M platform increases notablely along with the increase of the number of products and their order quantities. However, if the number of products and their order quantities are increased too much and far exceed the production capacity of each product of all manufacturers, the post-collaboration total profit will decrease. This is showed in Tables 8, where the order quantity of each product in each instance is 1.2 times of that in original instances set 5.

Table 8 Results of the instances set 5 (10x5) with larger order quantities

| **Instance** | **1** | **2** | **3** | **4** | **5** | **6** | **7** | **8** | **9** | **10** |
|---|---|---|---|---|---|---|---|---|---|---|
| Post-collaboration total profit | 518827 | 395625 | 500683 | 479818 | 465037 | 414284 | 476312 | 288983 | 435384 | 325876 |
| The value of $\gamma$ | 1 | 1 | 1 | 1 | 1 | 1 | 1 | 1 | 1 | 1 |
| Product shortage | Y | Y | Y | Y | Y | Y | Y | Y | Y | Y |

3) Our proposed mechanism can generate a profit allocation that stays in the core ($\gamma$ = 1) in 31 of 60 instances and yield a feasible allocation ($\gamma$-approximate core allocations) for all other 29 instances (0.73 < $\gamma$ < 1). The core allocations can guarantee the stability of the collaborative alliance and the sustainability of C2M platform. The $\gamma$-approximate core allocations can achieve the similar effect in some sense.

For each instance of instance set 2, Tables 9 list the post-collaboration profit of each manufacturer in detail. The profit of each manufacturer is increased through joining the collaboration of the C2M platform. With constraints (9) and (10), the manufacturer who fulfils more valuable production orders is allocated a larger profit after collaboration under the condition of core or $\gamma$-approximate core constraints. Thus, the fairness of our proposed mechanism can be demonstrated. If two manufacturers perform the same valuable production orders, they may be allocated the same profit.

Table 9 Post-collaboration profits of the manufacturers in instances set 2 (1x10)

| **Manufacturer** | **1** | **2** | **3** | **4** | **5** | **6** | **7** | **8** | **9** | **10** |
|---|---|---|---|---|---|---|---|---|---|---|
| instances 1 | 18942 | 18942 | 18942 | 20618 | 27595 | 27595 | 18942 | 18942 | 27595 | 20618 |
| instances 2 | 8103 | 7593 | 8436 | 7593 | 11321 | 11321 | 11321 | 14516 | 11321 | 8436 |
| instances 3 | 31348 | 26305 | 31348 | 26305 | 20794 | 28779 | 20794 | 28779 | 29126 | 29126 |
| instances 4 | 21729 | 17824 | 15155 | 15155 | 15155 | 21729 | 21729 | 24692 | 21729 | 24692 |
| instances 5 | 7602 | 7602 | 7602 | 7602 | 16062 | 7602 | 7602 | 7602 | 16062 | 7602 |
| instances 6 | 23746 | 30585 | 31675 | 17801 | 26768 | 17801 | 30585 | 24072 | 17801 | 30585 |
| instances 7 | 4817 | 3442 | 4262 | 2005 | 2005 | 1571 | 2005 | 3442 | 2452 | 4900 |
| instances 8 | 5939 | 6763 | 4547 | 5939 | 4547 | 6763 | 4059 | 6763 | 6763 | 5939 |
| instances 9 | 14061 | 8141 | 18466 | 14061 | 14061 | 21786 | 18466 | 8141 | 18349 | 14061 |
| instances 10 | 10930 | 19002 | 17240 | 24880 | 17240 | 24880 | 24880 | 10930 | 24880 | 10930 |

After a careful check on the data files and the experiment results, we also find that the manufacturer who has lower manufacturing cost and larger production capacity of each product are allocated more production order quantities and obtain more profit in collaboration.

The same observations can be found in other instance sets, they are not showed here due to the page limit.

## 5 Conclusion

Along with the changing of customers' consumption pattern, the customer-to-manufacturer (C2M) business model has attracted more and more attention. A profit allocation problem for C2M platform is studied in this paper. A $\gamma$-approximate core guaranteed profit allocation mechanism is proposed under a centralized collaboration framework, which can assure the stability of the alliance and the sustainability of C2M platform. Numerical experiments on sixty instances generated from data of Alibaba are tested to evaluate the mechanism and relevant models. Our future work is to explore other kind of collaboration frameworks and corresponding profit allocation models for C2M platform.

## Acknowledgment

This study is supported by the Social Science Achievements Appraisal Committee Foundation of Hunan Province (No. XSP21YBC477), the Scientific Research Project of Hunan Education Department (No. 19C1053), and the Innovation-Driven Foundation of Hunan University of Technology and Business (No. 2020QD02).